# TAIL ASYMPTOTICS FOR THE MAXIMUM OF PERTURBED RANDOM WALK

### By Victor F. Araman and Peter W. Glynn

#### *New York University and Stanford University*

Consider a random walk $S = (S_n : n \geq 0)$ that is "perturbed" by a stationary sequence $(\xi_n : n \geq 0)$ to produce the process $(S_n + \xi_n : n \geq 0)$. This paper is concerned with computing the distribution of the all-time maximum $M_\infty = \max\{S_k + \xi_k : k \geq 0\}$ of perturbed random walk with a negative drift. Such a maximum arises in several different applications settings, including production systems, communications networks and insurance risk. Our main results describe asymptotics for $\mathbb{P}(M_\infty > x)$ as $x \to \infty$. The tail asymptotics depend greatly on whether the $\xi_n$'s are light-tailed or heavy-tailed. In the light-tailed setting, the tail asymptotic is closely related to the Cramér–Lundberg asymptotic for standard random walk.

**1. Introduction.** Let $S = (S_n : n \geq 0)$ be a random walk sequence, so that $S_0 = 0$ and $S_n = X_1 + \cdots + X_n$ for $n \geq 1$, where $X_1, X_2, \ldots$ define an i.i.d. sequence of random variables. We set $\mathbb{E}X_1 = \mu$. We shall generally assume that the random walk has negative drift. Given a sequence $(\xi_n : n \geq 0)$ of "perturbations," we call the process $(S_n + \xi_n : n \geq 0)$ a "perturbed random walk." This paper is concerned with developing limit theorems and related approximations for $\mathbb{P}(M_\infty > x)$, where $M_\infty$ is the all-time maximum

$$M_\infty = \max\{S_k + \xi_k : k \geq 0\}.$$

Perturbed random walks have been previously studied within the insurance risk theory literature; see, for example, [7, 16, 17]. The existing literature typically assumes that the perturbation sequence $(\xi_n : n \geq 0)$ is itself of random walk type. However, in this paper, we shall instead assume that the perturbations $(\xi_n : n \geq 0)$ form a stationary sequence, so that the large-time behavior of the sequence $(S_n + \xi_n : n \geq 0)$ is basically inherited from that associated with $S$ itself. Such perturbed random walks (with stationary perturbations) arise naturally in several different applications settings. In Section 2, we argue that $M_\infty$ describes the steady-state order-to-delivery time in a make-to-order production system, in which delays are incurred at the production facility due to supplier tardiness in component delivery. The quantity $M_\infty$ also arises in the modeling of end-to-end delay of data packets in a telecommunications network, as well as in computing ruin probabilities for an insurance firm in which delayed premium payment occurs; see, for instance, [1].

Perturbed random walks bear some resemblance to the process studied within the nonlinear renewal theory literature; see [12, 13, 19] for a discussion of the main results. However, nonlinear renewal theory presumes that the perturbations are uniformly continuous in probability, in the sense that $\max\{|\xi_{n+k} - \xi_n| : 0 \leq k \leq n\delta\}$ converges to zero in probability as $n \to \infty$. This continuity condition is clearly violated by the i.i.d. perturbations that are a central example in our theory. Additional related literature includes [10], in which first passage times for perturbed random walk with positive drift are analyzed, as well as [8] and [9]. The latter two references obtain approximations for the tail of $M_\infty$ when the $\xi_n$'s are constant after a finite number of steps (and hence satisfy the uniform continuity in probability hypothesis).

The all-time maximum $M_\infty$ of our perturbed random walk inherits characteristics related to the all-time maximum $\max\{S_k : k \geq 0\}$ of ordinary random walk and from the extreme value behavior of the perturbations $(\xi_k : k \geq 0)$. Thus, the theory we develop here exhibits new features that are not present in the context of either nonlinear renewal theory or the perturbed random walks previously considered in the literature. In particular, it turns out that the tail behavior of the perturbations interacts with the (unperturbed) random walk in a nontrivial way to produce the tail behavior of $M_\infty$.

We establish the following results for the all-time maximum $M_\infty$:

1. We find an integral equation satisfied by the distribution of $M_\infty$ and characterize its solution (Proposition 3).
2. When the perturbations are light-tailed, we obtain exact asymptotics for $\mathbb{P}(M_\infty > x)$ as $x \to \infty$, under conditions closely related to those of the Cramér–Lundberg tail asymptotic for unperturbed random walk (Theorems 1 and 2).
3. For heavy-tailed perturbations, we obtain exact asymptotics for $\mathbb{P}(M_\infty > x)$ as $x \to \infty$ (Theorems 3 and 4).

We should note that under some regularity conditions, the results of Propositions 1, 2 and 6, as well as Theorem 4, hold not only for a random walk with i.i.d. increments, but also for a generalized random walk where the increments $X$ form a stationary sequence (see [1] for more details). Furthermore, in the presence of light-tailed perturbations and under weak conditions on the underlying (generalized) random walk, we can also prove an asymptotic for $\log \mathbb{P}(M_\infty > x)$ as $x \to \infty$ (we again refer the reader to [1] for more details).

Araman and Glynn [2] develop related results for the distribution of $M_\infty$ when the random walk $S$ has a negative drift close to zero, so that diffusion approximation techniques can be applied.

This paper is organized as follows. Section 2 focuses on a motivating example and basic properties of perturbed random walk, while Section 3 is concerned with developing asymptotics for $\mathbb{P}(M_\infty > x)$ when the $\xi_n$'s are light-tailed. Section 4 is devoted to heavy-tailed perturbations.

**2. A motivating example and basic properties.** Several models exist in which perturbed random walks arise naturally. We provide here a production model and direct the reader to [1] for a discussion of telecommunications and insurance risk examples.

EXAMPLE. Consider here a model of a so-called "make-to-order" production facility, which incorporates possible delays in the arrival of supplier components. Incorporating such delays is an important modeling issue and has been previously addressed by many authors, such as Kaplan [11], Sahin [15], Zipkin [20] and many others, using different modeling approaches than those we shall now describe.

For the purposes of our model, let $\tilde{A}_n$ be the time at which the $n$th order to the production facility is placed, so that $\tilde{A} = (\tilde{A}_n : n \geq 0)$ is nondecreasing. The facility satisfies demand according to a "first-come, first-served" priority rule, so that the $n$th order is processed to completion prior to initiating the $(n + 1)$st order, for $n \geq 0$. As soon as an order is placed, the production facility contacts its supplier. The supplier provides to the production facility the components necessary to produce the $n$th finished item after a delay of duration $\eta_n$, so that the corresponding components arrive to the production facility at time $A_n \stackrel{\triangle}{=} \tilde{A}_n + \eta_n$. Thus, we are modeling a "just-in-time" facility in which delays may be a consequence of transportation time, the time required to produce necessary components at the supplier facility, or some combination of these factors. Note that the sequence $A = (A_n : n \geq 0)$ need not be nondecreasing. Although the production facility satisfies customer demand according to the original order sequence, component deliveries to the producer may occur in an order that differs from this original ordering. This possibility differentiates our model from previous work on supply delay issues. Most of the references mentioned earlier model stochastic supply leadtimes as being either exogenous quantities (independent of the arrival process and the processing time) or endogenous in the sense that they are determined by an embedded queueing system. In both cases order-crossovers are not allowed.

As indicated above, customer demand is satisfied according to the original order sequence. This assumption is realistic in settings where customer equity concerns and total-order-to-delivery time issues predominate. Note, however, that, under this assumption, the production facility is no longer a "work-conserving" system (in the sense that the facility may be idle with orders present). If $D_n$ is the time at which the production facility completes work on the $n$th order, this assumption guarantees that the production facility initiates processing on the $(n+1)$st order at the maximum of $D_n$ and $A_{n+1}$. Let $V_{n+1}$ be the processing time required for the $(n+1)$st order, $\tilde{W}_{n+1}$ the time that an order spends waiting in the production facility and $W_n$ the total time spent in the system, so that $W_n = \tilde{W}_n + \eta_n$. A simple argument shows that

$$(2.1) \qquad \tilde{W}_{n+1} = \max(\tilde{W}_n + A_n - A_{n+1} + V_n, 0).$$

Thus, we recover the familiar Lindley recursion (2.1) (for single server queues), even though we are not requiring here that $A = (A_n : n \geq 0)$ be nondecreasing. Let $Z_{n+1} = \tilde{A}_n - \tilde{A}_{n+1} + V_n$ for $n \geq 0$ and assume, for simplicity, that the production facility is idle at time $t = 0$, and that the 0th order arrives at time $\tilde{A}_0 = 0$. Then (as is usual for the Lindley recursion), it is easily seen that

$$(2.2) \qquad \tilde{W}_n = \max\left\{ \sum_{j=k+1}^{n} Z_j + \eta_k - \eta_n : 0 \leq k \leq n \right\}.$$

Proposition 1 depicts the steady-state behavior of $W_n$ as $n \to \infty$.

CONDITION A1. $((Z_n, \eta_{n-1}) : n \geq 1)$ is a sequence of i.i.d. random vectors.

Let $((Z_n, \eta_{n-1}) : -\infty < n < +\infty)$ be a two-sided version of $((Z_n, \eta_{n-1}) : n \geq 1)$ and set $X_n = Z_{-n}$, $\xi_k = \eta_{-k-1}$ and $S_n = X_1 + X_2 + \cdots + X_n$ for $n \geq 1$ with $S_0 = 0$. It is then easy to see that

PROPOSITION 1.   *Under Condition* A1, $W_n \Rightarrow M_\infty$ *as* $n \to +\infty$, *where*

$$M_\infty = \max\{S_k + \xi_k : k \geq 0\}.$$

Because total order-to-delivery time is a key performance characteristic for a production facility, it follows that computing the distribution of $M_\infty$ is important, where

$$M_n = \max\{S_k + \xi_k : 0 \leq k \leq n\}.$$

We list below some basic results for which we have omitted the proofs; for the details, see [1]. Set $\xi_j^+ = \max(\xi_j, 0)$.

PROPOSITION 2. *Suppose that $(X_j : j \geq 1)$ is a sequence of i.i.d. random variables such that $\mathbb{E}X_1 = \mu$, where $-\infty < \mu < 0$.*

(i) *If $(\xi_j : j \geq 0)$ is an identically distributed sequence with $\mathbb{E}\xi_j^+ < \infty$, then $M_\infty < \infty$ a.s.*

(ii) *If $(\xi_j : j \geq 0)$ is a sequence of i.i.d. random variables, then $M_\infty < \infty$ a.s. implies that $\mathbb{E}\xi_j^+ < \infty$.*

To compute the exact distribution of $M_\infty$, define, for $g : \mathbb{R} \to \mathbb{R}^+$, the linear operator $T$ via

$$(Tg)(x) = \int_{\mathbb{R}} g(x - y) \mathbb{P}(\xi_0 \leq x, X_1 \in dy).$$

Also, let $T^0$ be the identity operator and define $T^n$ inductively via $T^{n+1} = T \circ T^n$ for $n \geq 0$.

PROPOSITION 3. *Suppose that $((\xi_j, X_{j+1}) : j \geq 0)$ is a sequence of i.i.d. random pairs.*

(i) *If $u^*(x) = \mathbb{P}(M_\infty > x)$, then $u^* = (u^*(x) : x \in \mathbb{R})$ satisfies the linear integral equation*

$$u = b + Tu,$$

*where $b = (b(x) : x \geq 0)$ is given by $b(x) = \mathbb{P}(\xi_0 > x)$.*

(ii) *Furthermore, $u^*$ is given by*

$$u^* = \sum_{n=0}^{\infty} T^n b. \tag{2.3}$$

The proof is similar in spirit to that for the unperturbed case discussed in [6]. There is one special case for which we have been able to derive a closed form expression for $u$.

PROPOSITION 4. *Suppose that $(-X_j : j \geq 1)$ is a sequence of i.i.d. exponential($\lambda$) random variables, independent of the i.i.d. perturbation sequence $(\xi_j : j \geq 0)$. Then,*

$$\mathbb{P}(M_\infty \leq x) = \mathbb{P}(\xi_0 \leq x) \exp\left(-\lambda \int_x^\infty \mathbb{P}(\xi_0 > y) \, dy\right)$$

*for $x \in \mathbb{R}$.*

Because we are generally unable to solve for the distribution of $M_\infty = \max\{S_n + \xi_n : n \geq 0\}$ exactly, we shall instead satisfy ourselves with obtaining approximations for the tail probability $\mathbb{P}(M_\infty > x)$. In [2], we obtain heavy-traffic diffusion approximations for this tail probability when the mean of the underlying random walk is negative but close to zero. In the current paper, we focus on studying the asymptotics of $\mathbb{P}(M_\infty > x)$ for $x$ large (when the mean of the random walk need not be close to zero).

**3. Tail asymptotics for light-tailed perturbations.** In this section we shall focus on the case in which the behavior of $\mathbb{P}(M_\infty > x)$ for large $x$ is primarily determined by the distribution of the random walk (as opposed to the perturbations). Specifically, we shall be interested in studying conditions on the perturbed random walk under which the tail of $M_\infty$ satisfies a Cramér–Lundberg asymptotic of the form

$$(3.1) \qquad\qquad \mathbb{P}(M_\infty > x) \sim c \exp(-\theta^* x)$$

as $x \to \infty$, where the rate constant $\theta^*$ agrees with that associated with the tail of the maximum of the unperturbed random walk as described through the conventional Cramér–Lundberg asymptotic; see, for example, [3] for a discussion of the relevant theory for unperturbed random walk. We refer to this setting as the "light-tailed" perturbation case.

Note that $\mathbb{P}(M_\infty > x) \geq \mathbb{P}(\xi_0 > x)$. Hence, if there exist positive constants $d$ and $\nu$ such that

$$(3.2) \qquad\qquad \mathbb{P}(\xi_j > x) \sim d \exp(-\nu x)$$

as $x \to \infty$, it follows that the validity of (3.1) requires that $\nu \geq \theta^*$. It follows that if the $\xi_j$'s satisfy (3.2) with $\nu \geq \theta^*$, then $\mathbb{E}\exp(\theta \xi_j) < \infty$ for $\theta < \theta^*$. In the remainder of this section, we make the (slightly) stronger hypothesis that $\mathbb{E}\exp(\theta \xi_j) < \infty$ for $\theta$ in a neighborhood of $\theta^*$.

CONDITION A2. $(X_j : j \geq 1)$ is an i.i.d. sequence of random variables for which there exist positive constants $\theta^*$ and $\varepsilon$ such that $\mathbb{E}\exp(\theta X_1) < \infty$ for $|\theta - \theta^*| < \varepsilon$ and $\mathbb{E}\exp(\theta^* X_1) = 1$.

It is well known that if the $X_i$'s are nonlattice, then Condition A2 ensures that the exact asymptotic

$$(3.3) \qquad\qquad \mathbb{P}\Big(\max_{n \geq 0} S_n > x\Big) \sim r \exp(-\theta^* x)$$

as $x \to \infty$ (for some positive $r$) holds for the maximum of the unperturbed random walk; see, for example, [3]. This, of course, is precisely the Cramér–Lundberg asymptotic for the (unperturbed) random walk $(S_n : n \geq 0)$.

A key idea in extending the exact asymptotic (3.3) to the context of perturbed random walk is the use of "change-of-measure" techniques. In particular, let $\mathbb{P}^*$ be the probability on the path space of the $(X_j, \xi_j)$'s under which the $X_j$'s are i.i.d. random variables with common distribution given by

$$\mathbb{P}^*(X_j \in dx) = \exp(\theta^* x)\mathbb{P}(X_j \in dx).$$

If $\mathbb{E}^*(\cdot)$ is the expectation operator corresponding to $\mathbb{P}^*$, and following a similar argument than in the unperturbed case (see [3]), we have that

$$(3.4) \qquad \exp(\theta^* x)\mathbb{P}(M_\infty > x) = \mathbb{E}^* \exp(-\theta^*(S_{T(x)} - x)),$$

where $T = T(x) \overset{\Delta}{=} \inf\{n \geq 0 : S_n + \xi_n > x\} < \infty$ $\mathbb{P}^*$-a.s.

Exact asymptotics for $\mathbb{P}(M_\infty > x)$ may therefore be deduced by studying the distribution of $S_{T(x)} - x$ under $\mathbb{P}^*$ for large values of $x$. In contrast to the corresponding analysis for unperturbed random walk, there is no guarantee here that $S_{T(x)} - x$ will be positive or that $T(x)$ will correspond to a ladder epoch of $(S_n : n \geq 0)$. As a consequence, the analysis of $S_{T(x)} - x$ under $\mathbb{P}^*$ is more intricate than in the setting of standard random walk. One approach to the analysis of $S_{T(x)} - x$ involves the use of renewal-theory ideas. For unperturbed random walk, the process $(S_{T(x)} - x : x \geq 0)$ is precisely the "current age" process corresponding to a renewal process in which the increments are governed by the strictly ascending ladder height of $(S_n : n \geq 0)$ under $\mathbb{P}^*$. However, in the setting of perturbed random walk, $(S_{T(x)} - x : x \geq 0)$ is not even a Markov process. Markov structure can be obtained by adding a supplementary state variable such as $\max_{0 \leq j \leq T(x)}(S_j + \xi_j) - x$. This two-dimensional Markov structure complicates the application of renewal theoretic ideas.

One setting in which one can directly apply such renewal-theoretic concepts is when the $X_j$'s are unbounded above while the $\xi_j$'s are upper bounded.

THEOREM 1. *Suppose that either:*

(i) *$((X_{j+1}, \xi_j) : j \geq 0)$ is a sequence of i.i.d. random vectors or*

(ii) *$(X_j : j \geq 1)$ is a sequence of i.i.d. random variables that is independent of the stationary sequence $(\xi_j : j \geq 0)$.*

*Assume, in either case, that Condition A2 holds and that $X_1$ has a nonlattice distribution. Suppose, in addition, that there exists $c_0 < \infty$ such that $\mathbb{P}(X_1 > c_0) > 0$ and $\mathbb{P}(\xi_0 \leq c_0) = 1$. Then, there exists a positive constant $c$ such that*

$$(3.5) \qquad \mathbb{P}(M_\infty > x) \sim c \exp(-\theta^* x)$$

*as $x \to \infty$.*

PROOF. Set $\tau_+(0) = 0$ and let $\tau_+(n) = \inf\{j > \tau_+(n-1) : S_j > S_{\tau_+(n-1)}\}$ be the $n$th strict ascending ladder height epoch of $(S_k : k \geq 0)$. In addition, set $\varsigma(0) = 0$ and let $\varsigma(n) = \inf\{j : S_j > S_k + c_0; j > k \geq \varsigma(n-1)\}$ be the $n$th strict ascending ladder epoch at which the ladder height increment is greater than $c_0$. Observe that $\varsigma(1) = \tau_+(\eta)$, where $\eta = \inf\{j \geq 1 : S_{\tau_+(j)} - S_{\tau_+(j-1)} > c_0\}$. The random variable $\eta$ is geometric, with mass function

$$\mathbb{P}^*(\eta = l) = \mathbb{P}^*(S_{\tau_+(1)} > c_0)\mathbb{P}^*(S_{\tau_+(1)} \leq c_0)^{l-1}.$$

But,

$$\mathbb{P}^*(S_{\tau_+(1)} > c_0) \geq \mathbb{P}^*(S_{\tau_+(1)} > c_0, \tau_+(1) = 1) = \mathbb{P}^*(X_1 > c_0) > 0,$$

and consequently $\mathbb{E}^*\eta < \infty$, so that $\varsigma(1) < \infty$ a.s. The key characteristic of $\varsigma(1)$ is that $T(S_{\varsigma(1)}) \geq \varsigma(1)$. To see this, note that $S_{\varsigma(1)} > S_j + c_0$ for $0 \leq j < \varsigma(1)$. Hence, $S_{\varsigma(1)} > S_j + \xi_j$ for $0 \leq j < \varsigma(1)$, so that $T(y) \geq \varsigma(1)$ for $y \geq S_{\varsigma(1)}$. Put

$$u(x) = \mathbb{E}^* \exp(-\theta^*(S_{T(x)} - x))$$

and

$$b(x) = \mathbb{E}^* \exp(-\theta^*(S_{T(x)} - x))I(S_{\varsigma(1)} > x).$$

Then,

$$u(x) = b(x) + \int_0^x \mathbb{E}^* \exp(-\theta^*(S_{T(x)} - x))I(S_{\varsigma(1)} \in dy)$$

$$= b(x) + \int_0^x \mathbb{E}^* \exp\left(-\theta^* \sum_{j=\varsigma(1)}^{\infty} (S_j - x)I(S_j + \xi_j > x,\right.$$

$$\left. S_k + \xi_k \leq x, k < j)\right)$$

$$\times I(S_{\varsigma(1)} \in dy)$$

$$= b(x) + \int_0^x \mathbb{E}^* \exp\left(-\theta^* \sum_{j=0}^{\infty}(S_{\varsigma(1)+j} - S_{\varsigma(1)} - (x-y))\right.$$

$$\times I(S_{\varsigma(1)+j} - S_{\varsigma(1)} + \xi_{\varsigma(1)+j} > x - y,$$

$$\left. S_{\varsigma(1)+k} + \xi_{\varsigma(1)+k} \leq x, k < j)\right)$$

$$\times I(S_{\varsigma(1)} \in dy)$$

$$= b(x) + \int_0^x u(x-y)F_\varsigma(dy),$$

where $F_\varsigma(dx) = \mathbb{P}^*(S_{\varsigma(1)} \le x)$. The fact that $\varsigma(1) \le T(x)$ on $\{S_{\varsigma(1)} \le x\}$ was used in the second equality above, whereas the independence between $(\xi_n, \xi_{n+1}, \ldots)$ and $S_n$ (together with the stationarity of the $\xi_j$'s) provides the third equality. [Note that the required independence and stationarity hold under either hypothesis (i) or (ii) of the theorem.]

It remains only to apply Feller's version of the renewal theorem to the above renewal equation. First, observe that because $X_1$ is nonlattice, $\mathbb{P}^*(S_{\tau_+(1)} \in \cdot)$ is nonlattice; see page 222 of [3]. For the integrability of $b(\cdot)$, we use the fact that $S_{T(x)} + \xi_{T(x)} > x$, so that $-\theta^*(S_{T(x)} - x) = -\theta^*(S_{T(x)} + \xi_{T(x)} - x) + \theta^*\xi_{T(x)} \le \theta^* c_0$. Consequently, Wald's identity implies that

$$\int_0^\infty b(x)\,dx \le e^{\theta^* c_0}\mathbb{E}^* S_{\varsigma(1)}$$
$$= e^{\theta^* c_0}\mathbb{E}^*\eta\,\mathbb{E}^* S_{\tau_+(1)} < \infty.$$

The final step is to prove that $b(\cdot)$ is directly Riemann integrable. The sample paths of $(S_{T(x)} - x : x \ge 0)$ are right continuous with left limits everywhere. Because $\exp(-\theta^*(S_{T(x)} - x)) \le \exp(\theta^* c_0)$, the bounded convergence theorem implies that $b(\cdot)$ is also right continuous with left limits. It follows that $b(\cdot)$ has at most countably many discontinuities; see page 116 of [5]. Owing to the bound $b(x) \le \exp(\theta^* c_0)\mathbb{P}^*(S_{\varsigma(1)} > x)$ and Proposition 5.4.1 of [3], we may therefore conclude that $b(\cdot)$ is directly integrable, proving the theorem. $\square$

Assumption (i) of Theorem 1 allows some dependency between the random walk and the perturbation. However, the conclusion of Theorem 1 generally does not hold if we assume only that the perturbation is light-tailed with $\mathbb{E}\exp(\theta\xi_j) < \infty$ for $\theta$ in a neighborhood of $\theta^*$; additional regularity, as presumed in Theorem 1, is needed. To see this, suppose that $X_i = R_i - \tilde{R}_i$ for $i \ge 1$, where $(R_i : i \ge 1)$ and $(\tilde{R}_i : i \ge 1)$ are independent sequences of i.i.d. exponentially distributed random variables, with parameters $\lambda_1$ and $\lambda_2$, respectively. Choose $\lambda_2 = \lambda_1/4$ and let $\theta^* = \lambda_1 - \lambda_2 = 3\lambda_1/4$. Suppose that the $\xi_i$'s are correlated with the random walk. In particular, if $\xi_i = R_i$ for $i \ge 1$, then $M_\infty \ge 2R_1 - \tilde{R}_1$. In this case $2R_1 - \tilde{R}_1$ has an exponential right tail with parameter $\lambda_1/2$. Consequently, $M_\infty$ has a tail that can decrease no more quickly than exponentially with rate constant $\lambda_1/2$, which is heavier than that predicted by Theorem 1.

When the $\xi_j$'s are independent of the unperturbed random walk [as in assumption (ii) of Theorem 1], and the random walk has spread-out increments, we can remove the boundedness assumption on the perturbations. To deal with the complications that arise here, we use a proof based on coupling ideas. Before stating the main result of this paper, we prove the following lemma related to the general theory of unperturbed random walk (that we have been unable to find in the existing literature):

Lemma 1. *Consider a random walk with an i.i.d. sequence of increments $(X_n : n \geq 0)$, in which $0 < \mathbb{E}X_1 < \infty$. The increment $X_1$ of the random walk is spread-out if and only if the increment of the ascending ladder height $S_{\tau_+(1)}$ is spread-out.*

Proof. Let $F$ be the distribution of $X_1$ and $F_+$ the distribution of $S_{\tau_+(1)}$. Let $N(x) = \inf\{n \geq 0 : S_n > x\}$. Assume first that $F$ is spread-out. Hence there exists $n > 0$ such that the distribution of $S_n$ admits a density component. For $z > 0$ and $\varepsilon > 0$, form the random variables $Y_1 = S_{\lfloor z/\varepsilon \rfloor} I(X_1 \leq -\varepsilon, X_2 \leq -\varepsilon, \ldots, X_{\lfloor z/\varepsilon \rfloor} \leq -\varepsilon)$ and $Y_2 = (S_{n+\lfloor z/\varepsilon \rfloor} - S_{\lfloor z/\varepsilon \rfloor}) \times I(\max_{0 \leq j \leq n}(S_{j+\lfloor z/\varepsilon \rfloor} - S_{\lfloor z/\varepsilon \rfloor}) \leq z)$. For sufficiently large $z$, the random variable $Y_2 \overset{D}{=} S_n I(\max_{0 \leq j \leq n} S_j \leq z)$ admits a density component and, hence, the sum $Y_1 + Y_2$ does as well. Finally, by writing

$$(3.6) \qquad \mathbb{P}(S_{\tau_+(1)} \leq b) \geq \int_{-\infty}^0 \mathbb{P}(Y_1 + Y_2 \in dx)\mathbb{P}(S_{N(-x)} \leq b)$$

for all $b > 0$, we conclude that $F_+$ admits a density component. For the converse, suppose that $F$ is not spread-out. Then, $F_+^{*n}$ is concentrated on a null Lebesgue set for all $n$. Since $\mathbb{P}(S_{\tau_+(1)} \in \cdot, \tau_+(1) = k) = \mathbb{P}(S_k \in \cdot, \tau_+(1) = k) \leq F^{*k}(\cdot)$, evidently $F_+$ is also concentrated on a null Lebesgue set. A similar argument establishes that $F_+^{*n}$ is also concentrated on a null Lebesgue set for all $n$. □

Our next theorem generalizes part (ii) of Theorem 1 [but not part (i)].

Theorem 2. *Assume that $(X_j : j \geq 1)$ satisfies Condition A2. Suppose that $(\xi_j : j \geq 0)$ is a stationary sequence independent of $(X_j : j \geq 1)$, for which $\mathbb{E}\exp(\theta\xi_j) < \infty$ for $\theta$ in a neighborhood of $\theta^*$. If $X_1$ is spread-out, then there exists a positive constant $c$ such that*

$$\mathbb{P}(M_\infty > x) \sim c\exp(-\theta^* x)$$

*as $x \to \infty$.*

Proof. In contrast to the proof of Theorem 1, our proof here proceeds via a coupling construction (see, e.g., [20]). We construct a doubly-infinite stationary version of $(S_n : n \geq 0)$, calling it $(S_n^* : -\infty < n < \infty)$, and random times $\sigma_1$ and $\sigma_2$, for which $S_{\sigma_1+k} = S_{\sigma_2+k}^*$ for $k \geq 0$. The process $(S_n^* : -\infty < n < \infty)$ is stationary in the sense that for each $x \in \mathbb{R}$,

$$(3.7) \qquad (S_{N^*(x)+j}^* - x : -\infty < j < \infty) \overset{D}{=} (S_j^* : -\infty < j < \infty),$$

where $N^*(x) = \inf\{n \geq 0 : S_n^* > x\}$. We then extend the stationary sequence $(\xi_j : j \geq 0)$ to a doubly-infinite stationary sequence and set $T^*(x) = \inf\{n : S_n^* +$ $\xi_n > x\}$. The proof concludes by showing that $T^*(x) = T(x)$, for $x$ sufficiently large, and that $S^*_{T^*(x)} - x \overset{D}{=} S^*_{T^*(0)}$. It follows that

$$S_{T(x)} - x \Rightarrow S^*_{T^*(0)}$$

as $x \to \infty$. In the presence of uniform integrability of $(\exp(-\theta^*(S_{T(x)} - x)) : x \geq 0)$ under $\mathbb{P}^*$, we then obtain, recalling equation (3.4), the conclusions of the theorem. Note that both $N^*(\cdot)$ and $T^*(\cdot)$ depend on the indexing of the points in the set $(S^*_n : -\infty < n < \infty)$. It is therefore essential that our construction of $(S^*_n : -\infty < n < \infty)$ produces point indexing that respects the stationarity relation (3.7) [rather than the weaker stationarity relation $(S^*_n : -\infty < n < \infty) \overset{D}{=} (S^*_n + x : -\infty < n < \infty)$].

We start by using the conventional coupling for finite mean renewal processes with spread-out increment distribution to couple the ascending ladder heights (which is enabled by Lemma 1). Specifically, we couple $(S_{\tau_+(k)} : k \geq 0)$ to a stationary version $(U^*_k : -\infty < k < \infty)$ of the strictly ascending ladder height sequence, thereby constructing a pair of random times $\eta_1$ and $\eta_2$ for which $S_{\tau_+(\eta_1 + k)} = U^*_{\eta_2 + k}$ for $k \geq 0$; see, for example, [14]. Set $S^*_0 = U^*_0$. Conditional on $(U^*_j : -\infty < j < \eta_2)$, condition in the ladder height epochs $(\tau^*_+(j) : j \leq \eta_2)$ so that $S^*_{\tau^*_+(j)} = U^*_j$ for $j \leq \eta_2$. Finally, given $(S^*_{\tau^*_+(j)}, \tau^*_+(j) : j \leq \eta_2)$, generate $(S^*_n : n \leq \tau^*_+(\eta_2))$ from the corresponding conditional distribution and set $S^*_{\sigma_1 + k} = S_{\sigma_2 + k}$ for $k \geq 0$, where $\sigma_1 = \tau_+(\eta_1)$ and $\sigma_2 = \tau^*_+(\eta_2)$. The above construction couples $(S_n : n \geq 0)$ to $(S^*_n : -\infty < n < \infty)$ and yields a process $(S^*_n : -\infty < n < \infty)$ that is stationary in the sense of (3.7); see [1] for additional details. Turning now to the analysis of $T^*(x)$, observe that $T^*(x) > -\infty$, since $(S^*_n + \xi^+_n)/n \to -\mathbb{E}^* X_1 < 0$ $\mathbb{P}^*$-a.s. as $n \to -\infty$, so that (at most) finitely many of the $(S^*_j + \xi^*_j : j < 0)$ can exceed $x$. Since $S^*_n + \xi^+_n \to +\infty$ $\mathbb{P}^*$-a.s. as $n \to +\infty$, it follows that $T^*(x)$ is finite-valued for all $x$. Furthermore, $S^*_{T^*(x)} - x \overset{D}{=} S^*_{T^*(0)}$. To see this, note that

$$S^*_{T^*(x)} - x$$
$$= \sum_{k=-\infty}^{\infty} (S^*_k - x) I(S^*_k + \xi_k > x, S^*_j + \xi_j \leq x, j < k)$$
$$= \sum_{l=-\infty}^{\infty} (S^*_{N^*(x)+l} - x)$$
$$\qquad \times I(S^*_{N^*(x)+l} + \xi_{N^*(x)+l} > x, S^*_{N^*(x)+j} + \xi_{N^*(x)+j} \leq x, j < l)$$
$$\overset{D}{=} \sum_{l=-\infty}^{\infty} (S^*_{N^*(x)+l} - x) I(S^*_{N^*(x)+l} + \xi_l > x, S^*_{N^*(x)+j} + \xi_j \leq x, j < l)$$

$$\overset{D}{=} \sum_{l=-\infty}^{\infty} S_l^* I(S_l^* + \xi_l > x, S_j^* + \xi_j \leq x, j < l)$$

$$= S_{T^*(0)}^*.$$

To obtain the first distributional equality, we use the stationarity and independence of the $\xi_j$'s from the $S_n$'s and $S_n^{**}$'s. The second distributional equality is a consequence of (3.7).

We now show that $S_{T(x)} - x = S_{T^*(x)}^* - x$ for $x$ sufficiently large. This follows easily from the fact that $T(x) \to +\infty$ and $T^*(x) \to +\infty$ $\mathbb{P}^*$-a.s. as $x \to \infty$. This implies that $S_{T(x)} = S_{T^*(x)}^*$ for sufficiently large $x$ so that $T(x) \geq \sigma_1$ and $T^*(x) \geq \sigma_2$. Finally, the required uniform integrability follows from Proposition 6 below.   □

We conclude this section with a couple of bounds for $\mathbb{P}(M_\infty > x)$. We note that both the lower bound and upper bound are multiples of $\exp(-\theta^* x)$ and consequently are within a constant factor of the correct value of $\mathbb{P}(M_\infty > x)$ uniformly in $x$.

PROPOSITION 5.   *Suppose that $(\xi_j : j \geq 0)$ is a stationary sequence independent of $(S_n : n \geq 0)$. Then,*

$$\mathbb{P}(M_\infty > x) \geq \int_{\mathbb{R}} \mathbb{P}\Big(\max_{n \geq 0} S_n > x - y\Big) \mathbb{P}(\xi_0 \in dy).$$

*If, in addition, (3.3) holds with $\mathbb{E} \exp(\theta^* \xi_0) < \infty$, then*

$$\liminf_{x \to \infty} e^{\theta^* x} \mathbb{P}(M_\infty > x) \geq r \mathbb{E} \exp(\theta^* \xi_0).$$

PROOF.   Let $J$ be the time at which the random walk attains its all-time maximum and note that both $J$ and $S_J$ are independent of $(\xi_n : n \geq 0)$. Clearly, $M_\infty \geq S_J + \xi_J$, so

$$\mathbb{P}(M_\infty > x) \geq \int_0^\infty \mathbb{P}(S_J + \xi_J > x | S_J = y) \mathbb{P}(S_J \in dy)$$

$$= \int_0^\infty \mathbb{P}(\xi_0 > x - y) \mathbb{P}(S_J \in dy)$$

$$= \int_{\mathbb{R}} \mathbb{P}(S_J > x - y) \mathbb{P}(\xi_0 \in dy).$$

For the second conclusion, note that

$$e^{\theta^* x} \int_{\mathbb{R}} \mathbb{P}(S_J > x - y) \mathbb{P}(\xi_0 \in dy)$$

$$= \int_{\mathbb{R}} e^{\theta^*(x-y)} \mathbb{P}\Big(\max_{n \geq 0} S_n > x - y\Big) e^{\theta^* y} \mathbb{P}(\xi_0 \in dy).$$

Applying (3.3) [and noting that it implies that $\exp(\theta^*(x-y))\mathbb{P}(\max_{n\geq 0} S_n > x-y)$ is uniformly bounded in $x$ and $y$], the hypothesis $\mathbb{E}\exp(\theta^* \xi_0) < \infty$ permits the dominated convergence theorem to be invoked, yielding the result. $\quad\square$

PROPOSITION 6. *Assume Condition* A2 *and let* $\kappa \in (0, \theta^*)$ *satisfy* $\psi'(\kappa) = 0$. *Suppose that* $(\xi_j : j \geq 0)$ *is a stationary sequence, independent of* $(S_n : n \geq 0)$, *for which* $\mathbb{E}\exp(\theta\xi_1) < \infty$ *in a neighborhood of* $\theta^*$. *Then, for* $\theta$ *in a neighborhood of* $\theta^*$,

$$\limsup_{x\to\infty} \mathbb{E}^* \exp(-\theta(S_{T(x)} - x))$$
$$\leq \frac{1}{\psi'(\theta^*)} \mathbb{E}\xi_1 \exp(\theta\xi_1) + \left(1 + \frac{\mathbb{E}\exp(\theta\xi_1)}{1 - \exp(\psi(\kappa))}\right).$$

PROOF. Set $v(x) = \mathbb{E}^* \exp(-\theta(S_{T(x)} - x))$ and note that

$$\begin{aligned}
v(x) &= \mathbb{E}^* \exp(-\theta(S_{T(x)} - x)) I(T(x) < \tau_+(1)) \\
&\quad + \mathbb{E}^* \exp(-\theta(S_{T(x)} - x)) I(S_{\tau_+(1)} > x, T(x) \geq \tau_+(1)) \\
&\quad + \mathbb{E}^* \exp(-\theta(S_{T(x)} - x)) I(S_{\tau_+(1)} \leq x, T(x) \geq \tau_+(1)).
\end{aligned} \tag{3.8}$$

Set $k(x) = k_1(x) + k_2(x)$, where

$$\begin{aligned}
k_1(x) &= \mathbb{E}^* \exp(-\theta(S_{T(x)} - x)) I(T(x) < \tau_+(1)) \\
k_2(x) &= \mathbb{E}^* \exp(-\theta(S_{T(x)} - x)) I(S_{\tau_+(1)} > x, T(x) \geq \tau_+(1)).
\end{aligned}$$

The third term on the right-hand side of (3.8) can be written as

$$\begin{aligned}
&\mathbb{E}^* \sum_{n=0}^{\infty} \exp(-\theta(S_{n+\tau_+(1)} - x)) \\
&\qquad \times I(S_{j+\tau_+(1)} + \xi_{j+\tau_+(1)} \leq x, 0 \leq j < n, S_{n+\tau_+(1)} + \xi_{n+\tau_+(1)} > x) \\
&\qquad \times I(T(x) \geq \tau_+(1), S_{\tau_+(1)} \leq x) \\
&\quad \leq \sum_{n=0}^{\infty} \int_{(0,x]} \mathbb{E}^* \exp(-\theta(S_{n+\tau_+(1)} - S_{\tau_+(1)} - (x-y))) \\
&\qquad\qquad \times I(S_{j+\tau_+(1)} - S_{\tau_+(1)} + \xi_{j+\tau_+(1)} \leq x - y, 0 \leq j < n, \\
&\qquad\qquad\qquad S_{n+\tau_+(1)} - S_{\tau_+(1)} + \xi_{n+\tau_+(1)} > x - y) \\
&\qquad\qquad \times I(S_{\tau_+(1)} \in dy) \\
&\quad = \int_{(0,x]} v(x-y) \mathbb{P}^*(S_{\tau_+(1)} \in dy),
\end{aligned}$$

where the independence and the stationarity of the perturbations and the random walk was used for the equality. Hence,

$$(3.9) \qquad v \le k + F_+ * v,$$

where $*$ denotes convolution and $F_+(dx) = \mathbb{P}^*(S_{\tau_+(1)} \in dx)$ for $x \ge 0$. Iterating (3.9) $n$ times, we find that

$$(3.10) \qquad v \le \sum_{j=0}^n F_+^{(j)} * k + F_+^{(n+1)} * v,$$

where $F_+^{(j)}$ denotes the $j$-fold convolution of $F_+$. But $(F_+^{(n)} * v)(x) = \mathbb{E}^* v(x - S_{\tau_+(n)})$, and $v(\cdot)$ is bounded on $[0, x]$. To verify the boundedness, note that, for $\beta > 0$,

$$\mathbb{E}^* \exp(-\theta S_{T(x)})$$

$$\le \sum_{n=0}^\infty \mathbb{E}^* \exp(-\theta S_n) I(S_n + \xi_n > x)$$

$$= \sum_{n=0}^\infty \int_{\mathbb{R}} \mathbb{E}^* [\exp(-\theta S_n); S_n > x - y] \mathbb{P}(\xi_1 \in dy)$$

$$= \sum_{n=0}^\infty \int_{\mathbb{R}} \mathbb{E} [\exp((\theta^* - \theta) S_n); \beta S_n > \beta(x - y)] \mathbb{P}(\xi_1 \in dy)$$

$$\le \sum_{n=0}^\infty \int_{\mathbb{R}} \mathbb{E} [\exp((\theta^* - \theta + \beta) S_n)] \exp(-\beta x + \beta y) \mathbb{P}(\xi_1 \in dy)$$

$$= \exp(-\beta x) \mathbb{E} \exp(\beta \xi_1) \sum_{n=0}^\infty (\mathbb{E} \exp((\theta^* - \theta + \beta) X_1)^n),$$

which can be made convergent by choosing $\theta > \beta > \theta - \theta^*$. Given the boundedness, it follows that $(F^{(n)} * v)(x) \to 0$ as $n \to \infty$ and, hence,

$$v \le U_+ * k,$$

where $U_+ = \sum_{j=0}^\infty F_+^{(j)}$ is the renewal kernel associated with the distribution $F_+$. To apply Smith's version of the renewal theorem (see, e.g., page 187 of [3]) we now need to verify that $k_1$ and $k_2$ are dominated by nonincreasing integrable functions. For $k_1$, observe that

$$k_1(x) \le \mathbb{E}^* \exp(-\theta(S_{T(x)} + \xi_{T(x)} - x) + \theta \xi_{T(x)}) I(T(x) < \tau_+(1))$$

$$\le \mathbb{E}^* \exp(\theta \xi_{T(x)}) I(T(x) < \tau_+(1))$$

$$\le \mathbb{E}^* \exp\left(\theta \max_{0 \le j < \tau_+(1)} \xi_j\right) I\left(\max_{0 \le j < \tau_+(1)} (S_j + \xi_j) > x\right)$$

$$\le \mathbb{E}^* \exp\Big(\theta \max_{0 \le j < \tau_+(1)} \xi_j\Big) I\Big(\max_{0 \le j < \tau_+(1)} \xi_j > x\Big),$$

where the last inequality uses the fact that $S_j \le 0$ for $j < \tau_+(1)$. Hence, $k_1$ is dominated by a nondecreasing function for which

$$\begin{aligned}
&\int_0^\infty \mathbb{E}^* \exp\Big(\theta \max_{0 \le j < \tau_+(1)} \xi_j\Big) I\Big(\max_{0 \le j < \tau_+(1)} \xi_j > x\Big)\, dx \\
&= \mathbb{E}^* \max_{0 \le j < \tau_+(1)} \xi_j \exp(\theta \xi_j) \\
&\le \mathbb{E}^* \sum_{j=0}^{\tau_+(1)-1} \xi_j \exp(\theta \xi_j) \\
&= \mathbb{E}^* \tau_+(1) \cdot \mathbb{E}^* \xi_1 \exp(\theta \xi_1) \\
&= \mathbb{E}^* \tau_+(1) \cdot \mathbb{E} \xi_1 \exp(\theta \xi_1) < \infty,
\end{aligned} \tag{3.11}$$

where the second last equality follows from Wald's identity, and the final equality utilizes the fact that the $\xi_j$'s have the same distribution under $\mathbb{P}$ as under $\mathbb{P}^*$. As for the function $k_2$, note that

$$\begin{aligned}
k_2(x) &= \mathbb{E}^* \exp(-\theta(S_{T(x)} - x)) I(S_{\tau_+(1)} > x, T(x) \ge \tau_+(1), S_{T(x)} > x) \\
&\quad + \mathbb{E}^* \exp(-\theta(S_{T(x)} - x)) I(S_{\tau_+(1)} > x, T(x) \ge \tau_+(1), S_{T(x)} \le x) \\
&\le \mathbb{P}^*(S_{\tau_+(1)} > x, T(x) \ge \tau_+(1), S_{T(x)} > x) \\
&\quad + \mathbb{E}^* \exp(\theta \xi_{T(x)}) I(S_{\tau_+(1)} > x, T(x) \ge \tau_+(1), S_{T(x)} \le x) \\
&\le \mathbb{P}^*(S_{\tau_+(1)} > x) \\
&\quad + \mathbb{E}^* \exp(\theta \xi_{T(x)}) I(S_{\tau_+(1)} > x, T(x) \ge \tau_+(1), S_{T(x)} \le x).
\end{aligned} \tag{3.12}$$

The second term on the right-hand side of (3.12) is dominated by

$$\begin{aligned}
&\sum_{n=0}^\infty \mathbb{E}^* \exp(\theta \xi_{\tau_+(1)+n}) \\
&\quad \times I(S_{\tau_+(1)} > x, T(x) = \tau_+(1) + n, S_{\tau_+(1)+n} - S_{\tau_+(1)} \le 0) \\
&\le \mathbb{E}^* \exp(\theta \xi_1) \sum_{n=0}^\infty \mathbb{E}^* I(S_{\tau_+(1)} > x, T(x) = \tau_+(1) + n, \\
&\hspace{8cm} S_{\tau_+(1)+n} - S_{\tau_+(1)} \le 0) \\
&\le \mathbb{E} \exp(\theta \xi_1) \mathbb{P}^*(S_{\tau_+(1)} > x) \sum_{n=0}^\infty \mathbb{P}^*(S_n \le 0) \\
&= \mathbb{E} \exp(\theta \xi_1) \mathbb{P}^*(S_{\tau_+(1)} > x) \sum_{n=0}^\infty \mathbb{E} \exp(\theta^* S_n) I(S_n \le 0)
\end{aligned}$$

$$\leq \mathbb{E}\exp(\theta\xi_1)\mathbb{P}^*(S_{\tau_+(1)} > x)\sum_{n=0}^{\infty}\mathbb{E}\exp(\kappa S_n)$$

$$= \mathbb{E}\exp(\theta\xi_1)\mathbb{P}^*(S_{\tau_+(1)} > x)(1-\exp(\psi(\kappa)))^{-1},$$

from which the result follows immediately via an invocation of Smith's renewal theorem. $\square$

REMARK 1.   If the $\xi_j$'s satisfy the conditions of Proposition 6 and are, in addition, nonnegative, then $T(x) = \tau_+(1)$ on $\{S_{\tau_+(1)} > x, T(x) \geq \tau_+(1)\}$. It follows, in this case, that $k_2(x) \leq \mathbb{E}^*\exp(-\theta(S_{\tau_+(1)} - x))I(S_{\tau_+(1)} > x)$, so the upper bound

(3.13)
$$\limsup_{x\to\infty}\mathbb{E}^*\exp(-\theta(S_{T(x)} - x))$$

$$\leq \frac{1}{\psi'(\theta^*)}\mathbb{E}\xi_1\exp(\theta\xi_1)$$

$$+ \frac{1}{\mathbb{E}^*S_{\tau_+(1)}}\mathbb{E}^*\int_0^{S_{\tau_+(1)}}\exp(-\theta(S_{\tau_+(1)} - x))\,dx$$

holds.

When $\theta = \theta^*$, we note that the second term on the right-hand side of (3.13) is then the Cramér–Lundberg constant $r$ for the unperturbed random walk. Hence, both our upper and lower bounds can be expressed directly in terms of $r$ and $\theta^*$.

**4. Tail asymptotics for heavy-tailed perturbations.**   In this section we develop a couple of tail asymptotics for $M_\infty$ in which the perturbations have a sufficiently heavy right tail that they largely govern the behavior of $M_\infty$'s tail. We first consider the setting in which the $\xi_j$'s have an exponential-like right tail, for which $\mathbb{E}\exp(\theta^*\xi_1) = \infty$.

THEOREM 3.   *Assume that $(\xi_j : j \geq 0)$ is a sequence of i.i.d. random variables satisfying (3.2). Suppose also that $(X_j : j \geq 1)$ is a sequence of i.i.d. random variables, independent of $(\xi_j : j \geq 0)$, for which $\mathbb{E}\exp(\nu X_1) < 1$. Then,*

$$\mathbb{P}(M_\infty > x) \sim d(1 - \mathbb{E}\exp(\nu X_1))^{-1}\exp(-\nu x)$$

*as $x \to \infty$.*

PROOF. We start by observing that $S_n \to -\infty$ a.s. by the strong law of large numbers. In addition, Markov's inequality yields

$$\mathbb{P}\left(\max_{n\geq 0} S_n > x\right) \leq \sum_{n=0}^{\infty} \mathbb{P}(S_n \geq x)$$

$$= e^{-\nu x} \mathbb{E} \sum_{n=0}^{\infty} \exp(\nu S_n) I(S_n > x).$$

But $\sum_{n=0}^{\infty} \exp(\nu S_n) I(S_n > x) \to 0$ a.s. as $x \to +\infty$, and $\sum_{n=0}^{\infty} \exp(\nu S_n)$ has expectation $(1 - \mathbb{E}\exp(\nu X_1))^{-1} < \infty$. So, the dominated convergence theorem proves that $\mathbb{P}(\max_{n\geq 0} S_n > x) = o(e^{-\nu x})$ as $x \to \infty$.

By conditioning on $(S_j : j \geq 0)$, we find that

$$\mathbb{P}(M_\infty \leq x) = \mathbb{E} \exp\left(\sum_{j=0}^{\infty} \log(1 - \bar{F}_\xi(x - S_j))\right).$$

Fix $\varepsilon > 0$ and choose $\Delta$ large enough so that

$$-(1+\varepsilon)\bar{F}_\xi(x) \leq \log(1 - \bar{F}_\xi(x)) \leq -(1-\varepsilon)\bar{F}_\xi(x)$$

for $x \geq \Delta$. Then, on $\{\max_{n\geq 0} S_n \leq x - \Delta\}$, $x - S_j \geq \Delta$ for all $j \geq 0$, so that

$$-(1+\varepsilon)\sum_{j=0}^{\infty} \bar{F}_\xi(x - S_j) \leq \sum_{j=0}^{\infty} \log(1 - \bar{F}_\xi(x - S_j))$$

(4.1)

$$\leq -(1-\varepsilon)\sum_{j=0}^{\infty} \bar{F}_\xi(x - S_j).$$

Observe that, because $\bar{F}_\xi(x) \leq d_0 \exp(-\nu x)$ for $x \geq 0$ (and some finite constant $d_0$),

$$\sum_{j=0}^{\infty} \bar{F}_\xi(x - S_j) \leq d_0 \exp(-\nu x) \sum_{j=0}^{\infty} \exp(\nu S_j).$$

But $\sum_{j=0}^{\infty} \exp(\nu S_j)$ is integrable, so $\sum_{j=0}^{\infty} \exp(\nu S_j) < \infty$ a.s. Hence, $\sum_{j=0}^{\infty} \bar{F}_\xi(x - S_j) = O(e^{-\nu x}) \to 0$ a.s. as $x \to \infty$. As a consequence,

$$e^{\nu x}\left(1 - \exp\left(-(1-\varepsilon)\sum_{j=0}^{\infty} \bar{F}_\xi(x - S_j)\right)\right)$$

$$= e^{\nu x}(1-\varepsilon)\sum_{j=0}^{\infty} \bar{F}_\xi(x - S_j) + o(1) \qquad \text{a.s.}$$

$$\to (1-\varepsilon)d\sum_{j=0}^{\infty} \exp(\nu S_j) \qquad\qquad \text{a.s.}$$

as $x \to \infty$. A similar argument holds for the lower bound in (4.1). Since $\varepsilon > 0$ was arbitrary, this proves that

$$e^{\nu x}\left(1 - \exp\left(\sum_{j=0}^{\infty} \log(1 - \bar{F}_\xi(x - S_j))\right)\right) \to d\sum_{j=0}^{\infty} \exp(\nu S_j) \qquad \text{a.s.}$$

as $x \to \infty$.

To complete the proof, we need to show that we can invoke the dominated convergence theorem. This follows from

$$e^{\nu x}\left(1 - \exp\left(\sum_{j=0}^{\infty} \log(1 - \bar{F}_\xi(x - S_j))\right)\right)$$

$$\leq e^{\nu x} I(M > x - \Delta) + e^{\nu x}\left(1 - \exp\left(-(1 + \varepsilon)\sum_{j=0}^{\infty} \bar{F}_\xi(x - S_j)\right)\right)$$

$$\leq e^{\nu x} \cdot e^{-\nu(x - \Delta)} e^{\nu M} + e^{\nu x}(1 + \varepsilon)\sum_{j=0}^{\infty} \bar{F}_\xi(x - S_j)$$

$$\leq e^{\nu \Delta}\sum_{j=0}^{\infty} \exp(\nu S_j) + (1 + \varepsilon)d_0 \sum_{j=0}^{\infty} \exp(\nu S_j),$$

which, as argued earlier, has finite expectation. [The elementary bound $\exp(-y) \geq 1 - y$, for $y \geq 0$, was used to obtain the second inequality above.] □

Note that the constant $(1 - \mathbb{E}\exp(\nu X_1))^{-1}$ depends on the distribution of the increments of the random walk through an exponential moment. Our second "heavy-tailed perturbation" asymptotic shows that when the perturbations have a heavier tail than that described in Theorem 6, the distribution of $X_1$ affects the asymptotic only through its (ordinary) mean.

CONDITION A3.   $(X_j : j \geq 1)$ is an i.i.d. sequence of random variables for which there exists $\varepsilon > 0$ such that $\mathbb{E}\exp(\varepsilon|X_1|) < \infty$ and $\mathbb{E}X_1 < 0$.

In our next theorem, we assume that the perturbations have a hazard rate that decreases to zero. This assumption is one way of describing "heavy tails." In particular, all Pareto random variables, as well as all Weibulls random variables with shape parameter strictly less than one, obey this hazard rate hypothesis.

THEOREM 4.   *Assume Condition* A3. *Suppose that* $(\xi_j : j \geq 0)$ *is a sequence, independent of* $(S_j : j \geq 0)$, *of i.i.d. random variables having a continuous hazard rate function* $h(\cdot)$ *such that* $h(x) \to 0$ *as* $x \to \infty$. *Then,*

$$\mathbb{P}(M_\infty > x) \sim \frac{1}{|\mathbb{E}X_1|} \int_x^\infty \mathbb{P}(\xi_1 > y)\, dy$$

*as* $x \to \infty$.

PROOF. Set $R(x) = \int_x^\infty \bar{F}_\xi(y)\,dy$. Recalling that $h(x) \to 0$ and

$$(4.2) \qquad \frac{\bar{F}_\xi(x+y)}{\bar{F}_\xi(x)} = \exp\left(-\int_0^y h(x+u)\,du\right) \to 1$$

as $x \to \infty$, one can easily prove that

$$(4.3) \qquad R(x)/R(x+y) \to 1$$

as $x \to \infty$. Recall that $S_n/n \to \mu$ a.s. as $n \to \infty$ and observe that Condition A3 implies that, for $\varepsilon > 0$ and sufficiently small, there exists $\delta > 0$ such that

$$(4.4) \qquad \mathbb{P}(|S_n - n\mu| > \varepsilon n) = O(e^{-\delta n});$$

see page 18 of [4], for example. The Borel–Cantelli lemma therefore shows that $|S_n - n\mu| > \varepsilon n$ almost surely occurs only finitely often.

Fix $\varepsilon > 0$ and note that $L(\varepsilon) = \sup\{n \geq 0 : |S_n - n\mu| > \varepsilon n\}$ is therefore finite-valued. We may thus write

$$
\begin{aligned}
\sum_{j=0}^\infty \bar{F}_\xi(x - S_j) &= \sum_{j=0}^{L(\varepsilon)} \bar{F}_\xi(x - S_j) + \sum_{j=L(\varepsilon)+1}^\infty \bar{F}_\xi(x - S_j) \\
(4.5) \qquad &\leq (L(\varepsilon) + 1)\bar{F}_\xi\left(x - \max_{n \geq 0} S_n\right) + \sum_{j=L(\varepsilon)+1}^\infty \bar{F}_\xi(x - \mu j - \varepsilon j) \\
&\leq (L(\varepsilon) + 1)\bar{F}_\xi\left(x - \max_{n \geq 0} S_n\right) \\
&\quad + \frac{1}{|\mu + \varepsilon|} \sum_{j=L(\varepsilon)+1}^\infty \int_{x-(\mu+\varepsilon)(j-1)}^{x-(\mu+\varepsilon)j} \bar{F}_\xi(y)\,dy \\
&= (L(\varepsilon) + 1)\bar{F}_\xi\left(x - \max_{n \geq 0} S_n\right) + \frac{1}{|\mu + \varepsilon|} R(x - (\mu + \varepsilon)L(\varepsilon)) \\
&= \frac{1}{|\mu + \varepsilon|} R(x) + o(R(x)) \qquad \text{a.s.}
\end{aligned}
$$

as $x \to \infty$, where (4.2) to (4.4) are exploited for the final equality above. Since

$$(4.6) \qquad \frac{1}{R(x)}\mathbb{P}(M_\infty > x) = \mathbb{E}\frac{1}{R(x)}\left(1 - \exp\sum_{j=0}^\infty \ln(1 - \bar{F}_\xi(x - S_j))\right),$$

we focus on the random variable appearing on the right-hand side of (4.6). Equation (4.1) and (4.5) permit us to conclude that

$$\limsup_{x \to \infty} \frac{1}{R(x)} \left( 1 - \exp \sum_{j=0}^{\infty} \ln(1 - \bar{F}_\xi(x - S_j)) \right) \leq \frac{1}{|\mu + \varepsilon|} \qquad \text{a.s.}$$

A similar lower bound for the limit infimum holds. Since $\varepsilon > 0$ can be chosen arbitrarily small, we find that

$$(4.7) \qquad \frac{1}{R(x)} \left( 1 - \exp \sum_{j=0}^{\infty} \ln(1 - \bar{F}_\xi(x - S_j)) \right) \to \frac{1}{|\mu|} \qquad \text{a.s.}$$

as $x \to \infty$.

To complete the proof, it remains only to show that the left-hand side of (4.7) is uniformly integrable. Note that

$$\frac{1}{R(x)} \left( 1 - \exp \sum_{j=0}^{\infty} \ln(1 - \bar{F}_\xi(x - S_j)) \right) = \frac{1}{R(x)} \left( 1 - \prod_{j=0}^{\infty} (1 - \bar{F}_\xi(x - S_j)) \right).$$

For $0 \leq a_i \leq 1$, $i \geq 0$, an easy induction on $n$ shows that $1 - \prod_{i=0}^n (1 - a_i) \leq \sum_{i=0}^n a_i$ and, hence, the left-hand side is dominated by

$$
\begin{aligned}
&\leq \frac{1}{R(x)} I\left( \max_{n \geq 0} S_n > x/2 \right) \\
&\quad + \frac{1}{R(x)} \Bigg[ (L(\varepsilon_0) + 1) \bar{F}_\xi \left( x - \max_{n \geq 0} S_n \right) \\
&\qquad\qquad + \frac{R(x - (\mu + \varepsilon_0) L(\varepsilon_0))}{|\mu + \varepsilon_0|} \Bigg] \cdot I\left( \max_{n \geq 0} S_n \leq x/2 \right) \\
&\leq \frac{1}{R(x)} \exp\left( -\varepsilon_0 x/2 + \varepsilon_0 \max_{n \geq 0} S_n \right) I\left( \max_{n \geq 0} S_n > x/2 \right) \\
&\quad + \frac{1}{R(x)} (L(\varepsilon_0) + 1) \bar{F}_\xi \left( x - \max_{n \geq 0} S_n \right) I\left( \max_{n \geq 0} S_n > x/2 \right) + 1/|\mu + \varepsilon_0| \\
&\leq \frac{\exp(-\varepsilon_0 x/2)}{R(x)} \exp\left( \varepsilon_0 \max_{n \geq 0} S_n \right) + 1/|\mu + \varepsilon_0| \\
&\quad + (L(\varepsilon_0) + 1) \exp\left( \int_{x - \max_{n \geq 0} S_n}^{x} h(u)\, du \right) \cdot I\left( \max_{n \geq 0} S_n \leq x/2 \right) \frac{\bar{F}_\xi(x)}{R(x)} \\
&\leq \frac{\exp(-\varepsilon_0 x/2)}{R(x)} \exp\left( \varepsilon_0 \max_{n \geq 0} S_n \right) + 1/|\mu + \varepsilon_0| \\
&\quad + (L(\varepsilon_0) + 1) \exp\left( \sup\{h(y) : y \geq x/2\} \cdot \max_{n \geq 0} S_n \right) \cdot \max_{y \geq 0} \frac{\bar{F}_\xi(y)}{R(y)},
\end{aligned}
$$

(4.8)

for $\varepsilon_0 > 0$. But $\exp(-\varepsilon_0 x/2)/R(x) = \exp(-\int_0^x (\varepsilon_0/2 - h(u))\,du) \to 0$ as $x \to \infty$. It follows from (4.8) that the required uniform integrability of the left-hand side of (4.7) is a consequence of the Cauchy–Schwarz inequality and that

$$(4.9) \qquad \mathbb{E}\exp\left(2\varepsilon_0 \max_{n \geq 0} S_n\right) < \infty \quad \text{and} \quad \mathbb{E}L(\varepsilon_0)^4 < \infty$$

for some $\varepsilon_0 > 0$. But $\exp(2\varepsilon\max_{n\geq 0} S_n) \leq \sum_{n=0}^{\infty}\exp(\varepsilon S_n)$. Since $\psi'(0) < 0$, there exists $\theta > 0$ such that $\psi(\theta) < 0$. Condition A3 then ensures that $\mathbb{E}\sum_{n=0}^{\infty}\exp(\theta S_n) < \infty$. To deal with $L(\varepsilon_0)$, note that

$$L(\varepsilon)^4 \leq \sum_{j=0}^{L(\varepsilon)} j^4 I(|S_j - \mu j| > \varepsilon j) \leq \sum_{j=0}^{\infty} j^4 I(|S_j - \mu j| > \varepsilon j).$$

But (4.4) ensures that $\mathbb{P}(|S_j - \mu j| > \varepsilon j)$ decays exponentially in $n$, proving the finiteness of $\mathbb{E}L(\varepsilon)^4$. By choosing $\varepsilon_0$ sufficiently small so that $0 < \varepsilon_0 < \theta$ and $\mathbb{E}L(\varepsilon_0)^4 < \infty$, we obtain (4.9), completing the proof. $\square$

STERN SCHOOL OF BUSINESS
NEW YORK UNIVERSITY
44 WEST FOURTH STREET
NEW YORK, NEW YORK 10012
USA
E-MAIL: varaman@stern.nyu.edu

MANAGEMENT SCIENCE AND ENGINEERING
STANFORD UNIVERSITY
TERMAN ENGINEERING BUILDING
STANFORD, CALIFORNIA 94305
USA
E-MAIL: glynn@stanford.edu